\documentclass[12pt]{article}

\usepackage{pstricks,pst-node}
\usepackage{psfig}
\usepackage{epsfig}
\usepackage{amsmath,amsthm}
\usepackage{amssymb}
\usepackage{color}
\usepackage{xspace}
\usepackage[colorlinks=true,
linkcolor=green,
filecolor=brown,
citecolor=green]{hyperref}

\def\v{\vert}

\def\si{\sigma}

\newcommand{\seqnum}[1]{\href{http://oeis.org/#1}{\underline{#1}}}

\begin{document}
\newtheorem{theorem}{Theorem}
\newtheorem{defn}[theorem]{Definition}
\newtheorem{lemma}[theorem]{Lemma}
\newtheorem{prop}[theorem]{Proposition}
\newtheorem{cor}[theorem]{Corollary}
%\vspace*{-5mm}
\begin{center}
{\Large
The Number of $\bar{2}413\bar{5}$-Avoiding Permutations    \\ 
}

\vspace*{5mm}

DAVID CALLAN  \\
Department of Statistics  \\
\vspace*{-2mm}
University of Wisconsin-Madison  \\
\vspace*{-2mm}
1300 University Ave  \\
\vspace*{-2mm}
Madison, WI \ 53706-1532  \\
{\bf callan@stat.wisc.edu}  \\
\vspace*{5mm}
\end{center}

\begin{abstract}
We answer a question of R.\,J. Mathar and confirm that the counting 
sequence for $\overline{2}413\overline{5}$-avoiding permutations is 
the Invert transform of the Bell numbers. The proof relies on a 
simple decomposition of these permutations and the known fact that  
$\overline{2}413$-avoiding permutations are counted by the Bell 
numbers.

\end{abstract}

\section{Introduction} \vspace*{-5mm}
A permutation $\pi$ avoids the barred pattern $\overline{2}413\overline{5}$ 
if each instance of a not-necessarily-consecutive 413 pattern in $\pi$ 
is part of a 24135 pattern in $\pi$, and similarly for other barred patterns. 
Lara Pudwell \cite{schemes} presents a general approach to counting 
permutations avoiding a given 5-letter barred pattern that often 
produces a recurrence relation but does not do so in this case.
R.\,J. Mathar observed \cite{mathar} that the first few terms of 
the counting sequence for $\overline{2}413\overline{5}$-avoiding permutations 
agree with those of the Invert transform of the Bell 
numbers---Invert$(1,1,2,5,15,52,\ldots)=(1, 2, 5, 14, 43, 
144,\ldots)$---and asked if the two sequences 
coincide. We will show that the answer is yes. Section 2 reviews
terminology. Section 3 presents a decomposition for 
$\overline{2}413\overline{5}$-avoiding permutations in terms of 
$\overline{2}413$-avoiding permutations, yielding a 
bijection that proves the result.

\section{Review of terminology}  \vspace*{-5mm}
The Invert transform of a sequence $(a_{n})_{n\ge 1}$ is 
$(b_{n})_{n\ge 1}$ defined by
\[
1+\sum_{n\ge 1}b_{n}x^{n} = \frac{1}{1-\sum_{n\ge 1}a_{n}x^{n}},
\]
and has the following combinatorial interpretation \cite{canonical,cameron}.
If the counting sequence by size of a class of combinatorial 
structures, say A-structures, is $(a_{n})_{n\ge 1}$, then $b_{n}$ is 
the number of lists (of unspecified length) of A-structures whose 
total size is $n$.

For any barred pattern $\rho$, we use $S_{n}(\rho$) for the set of 
$\rho$-avoiding permutations of $[n]$. 
A permutation is \emph{standard} if its support set is an initial 
segment of the positive integers (or empty). To \emph{standardize} a 
permutation means to replace its smallest entry by 1, next smallest 
by 2, and so on. We use stand($\pi$) for the result of 
standardizing $\pi$.

\section{A decomposition and bijection}  \vspace*{-5mm}
We begin with two observations about a
$\overline{2}413\overline{5}$-avoiding permutation $\pi$. The entries 
\emph{after} $n$ in $\pi$ must decrease, else $n$ would start a 413 pattern 
with no available ``5''. If entries $c>a$ occur in that order \emph{before} 
$n$, then all elements of the interval $[a,c]$ must occur before $n$, 
else an element $b$ of $(a,c)$ would occur after $n$ and $cab$ is a 
413 pattern, again with no available ``5''.
From these observations, it follows that $\pi $ has the form
\[
\tau_{1}\tau_{2}\ \ldots \tau_{r}\,n\,a_{r-1}\,a_{r-2}\, \ldots \,a_{1}
\]
where

\vspace*{-5mm}
\begin{itemize}
    \item  each $\tau_{i}$ is a subpermutation, possibly empty, with support an interval of 
integers,
\vspace*{-2mm}
    \item  each $a_{i}$ is a single entry in $\pi$ and is 
    the only integer lying (in value) between the support intervals 
    of $\tau_{i}$ and $\tau_{i+1}, \ 1\le i \le r-1$,
\vspace*{-2mm}
    \item  $\tau_{1}<\tau_{2}<\ldots <\tau_{r}$ in the sense that each entry of 
$\tau_{i}$ is less than each entry of $\tau_{i+1},\ i=1,2,\ldots,r-1$,
\vspace*{-2mm}
    \item  $a_{r-1}>a_{r-2}> \ldots >a_{1}$,
\vspace*{-2mm}
    \item  each $\tau_{i}$ is $\overline{2}413$-avoiding.
\end{itemize}
\vspace*{-4mm} 
Conversely, any permutation $\pi$ with a decomposition $\tau_{1}\tau_{2}\ \ldots \tau_{r}\,n\,a_{r-1}\,a_{r-2}\, \ldots \,a_{1}$ satisfying these conditions is $\overline{2}413\overline{5}$-avoiding.

It is helpful to interpret the bulleted conditions in the setting of a modified 
permutation matrix where, for a permutation $\pi$, the entry in the 
$(i,j)$ cell, measuring from the southwest corner, is $\pi(i)$ if 
$j=\pi(i)$ and 0 otherwise. We see that a permutation $\pi$ has a decomposition that meets all the bulleted 
conditions if and only if its matrix has the form pictured 
schematically below for $r=4$ with each $\tau_{i}$ a $\overline{2}413$ 
avoider, each $a_{i}$ an entry of $\pi$, and 0's in all unshaded 
regions.

\begin{center}

\begin{pspicture}(-2,-0.7)(10,7.5)%\showgrid
\psset{unit=.6cm} 
\newgray{verylightgray}{.9} 

\pspolygon[linecolor=black,fillstyle=solid,fillcolor=white](0,0)(0,12)(12,12)(12,0)
\psline(2,0)(12,0)(12,2)(0,2)(0,3)(12,3)(12,5)(0,5)(0,6)(12,6)(12,8)(0,8)(0,9)(12,9)(12,11)(0,11)
\psline(2,0)(2,12)(4,12)(4,0)(6,0)(6,12)(8,12)(8,0)(9,0)(9,12)(10,12)(10,0)(11,0)(11,12)

\pspolygon[linecolor=black,fillstyle=solid,fillcolor=verylightgray](0,0)(0,2)(2,2)(2,0)
\pspolygon[linecolor=black,fillstyle=solid,fillcolor=verylightgray](2,3)(4,3)(4,5)(2,5)
\pspolygon[linecolor=black,fillstyle=solid,fillcolor=verylightgray](4,6)(6,6)(6,8)(4,8)
\pspolygon[linecolor=black,fillstyle=solid,fillcolor=verylightgray](6,9)(8,9)(8,11)(6,11)

\pspolygon[linecolor=black,fillstyle=solid,fillcolor=verylightgray](8,11)(9,11)(9,12)(8,12)
\pspolygon[linecolor=black,fillstyle=solid,fillcolor=verylightgray](9,8)(10,8)(10,9)(9,9)
\pspolygon[linecolor=black,fillstyle=solid,fillcolor=verylightgray](10,5)(11,5)(11,6)(10,6)
\pspolygon[linecolor=black,fillstyle=solid,fillcolor=verylightgray](11,2)(12,2)(12,3)(11,3)

\rput(1,1){\textrm{{\small $\tau_{1}$}}}
\rput(3,4){\textrm{{\small $\tau_{2}$}}}
\rput(5,7){\textrm{{\small $\tau_{3}$}}}
\rput(7,10){\textrm{{\small $\tau_{4}$}}}
\rput(8.5,11.5){\textrm{{\small $n$}}}
\rput(9.5,8.5){\textrm{{\small $a_{3}$}}}
\rput(10.5,5.5){\textrm{{\small $a_{2}$}}}
\rput(11.5,2.5){\textrm{{\small $a_{1}$}}}

\rput(6,-0.8){\textrm{{\small a $\overline{2}413\overline{5}$ avoider 
as a permutation matrix}}}

\end{pspicture}
\end{center}

With $a_{r}:=n$, the map $\pi \rightarrow 
\big(\textrm{stand}(\tau_{1}a_{1}),\ldots,\textrm{stand}(\tau_{r}a_{r})\big)$ 
is a bijection from $S_{n}(\overline{2}413\overline{5})$ 
to lists $(\si_{1},\ldots,\si_{r})$ with $r\ge 1$ where each $\si_{i}$ 
is an \emph{end-max avoider}---a standard $\overline{2}413\overline{5}$-avoiding 
permutation that ends at its maximum entry---and the total length of 
the $\si_{i}$'s is $n$. Clearly, the number of end-max avoiders of 
length $k$ is  the number of $\overline{2}413$-avoiding permutations 
of length $k-1$ and it is known \cite{eigensequence} that 
$\v\,S_{k-1}(\overline{2}413)\,\v =B_{k-1}$, the Bell number.

Set $a_{n}=B_{n-1}$, so that $a_{n}$ is the number of end-max avoiders of 
length (size) $n$. Then the Invert transform $(b_{n})_{n\ge 1}$ of $(a_{n})_{n\ge 1}$ is 
the number of lists of of end-max avoiders of total size $n$, which 
the bijection above shows is $\v\,S_{n}(\overline{2}413\overline{5})\,\v$. 
Hence, the counting sequence for $S_{n}(\overline{2}413\overline{5})$ 
is the Invert transform of $(a_{n})_{n\ge 1}=(1,1,2,5,15,\ldots)$, the 
full sequence of Bell numbers.


\begin{thebibliography}{99}
    
\bibitem{schemes}  Lara Pudwell, Enumeration Schemes for Permutations Avoiding 
Barred Patterns, 
\htmladdnormallink{\emph{Electronic J. Combinatorics}}{http://www.combinatorics.org/},
 \textbf{17} (1) (2010), R29.
  
\bibitem{mathar} R.\,J. Mathar, Comment on sequence \seqnum{A160701} in the 
\htmladdnormallink{On-Line Encyclopedia}{http://oeis.org/Seis.html} of Integer Sequences.   

\bibitem{canonical}  M. Bernstein and N.\,J.\,A. Sloane,  Some Canonical 
Sequences of Integers,
\emph{Linear Algebra Appl.} \textbf{226--228} (1995) 57--72. 

\bibitem{cameron} P.\,J. Cameron, Some sequences of integers, \emph{Discrete Math.}, \textbf{75} (1989) 89Ð102.

\bibitem{eigensequence} David Callan, A combinatorial interpretation of the eigensequence for composition,
\htmladdnormallink{\emph{J. Integer Sequences}}{http://www.cs.uwaterloo.ca/journals/JIS/}, \textbf{9}, Article 06.1.4, 2006. 

\end{thebibliography}
\end{document}